\documentclass{amsart}
\usepackage[T1]{fontenc}

\usepackage{amsmath,amsfonts,amssymb,amscd,amsthm,amsbsy,upref,mathrsfs}

\newtheorem{theorem}{Theorem}[section]
\newtheorem{lemma}[theorem]{Lemma}
\newtheorem{proposition}[theorem]{Proposition}
\newtheorem{corollary}[theorem]{Corollary}

\newtheorem{fact}[theorem]{Fact}
\newtheorem*{claim}{Claim}

\theoremstyle{definition}
\newtheorem{definition}[theorem]{Definition}
\newtheorem{notation}[theorem]{Notation}

\theoremstyle{remark}
\newtheorem{remark}[theorem]{Remark}

\newcommand{\bl}{\begin{lemma}}
\newcommand{\el}{\end{lemma}}
\newcommand{\bfa}{\begin{fact}}
\newcommand{\efa}{\end{fact}}
\newcommand{\bpr}{\begin{proposition}}
\newcommand{\epr}{\end{proposition}}
\newcommand{\bp}{\begin{proof}}
\newcommand{\ep}{\end{proof}}
\newcommand{\bd}{\begin{definition}}
\newcommand{\ed}{\end{definition}}
\newcommand{\bt}{\begin{theorem}}
\newcommand{\et}{\end{theorem}}
\newcommand{\bc}{\begin{corollary}}
\newcommand{\ec}{\end{corollary}}
\newcommand{\bn}{\begin{notation}}
\newcommand{\en}{\end{notation}}
\newcommand{\br}{\begin{remark}}
\newcommand{\er}{\end{remark}}
\newcommand{\bcl}{\begin{claim}}
\newcommand{\ecl}{\end{claim}}

\newcommand{\N}{{\mathbb{N}}}

\newcommand{\Q}{{\mathbb{Q}}}

\newcommand{\nrm}[1]{\|#1\|}
\newcommand{\al}{\alpha}
\newcommand{\e}{\varepsilon}
\newcommand{\de}{\delta}
\newcommand{\bnum}{\begin{enumerate}}
\newcommand{\enum}{\end{enumerate}}
\newcommand{\mc}{\mathcal}
\newcommand{\mt}{\mc{T}}
\newcommand{\mf}{\mc{F}}
\newcommand{\mg}{\mc{G}}
\newcommand{\ms}{\mc{S}}
\newcommand{\ma}{\mc{A}}
\newcommand{\mr}{\mc{R}}

\numberwithin{subsection}{section}
\numberwithin{equation}{section}

\newcommand{\norm}[1][\cdot]{\lVert #1\rVert}
\DeclareMathOperator{\supp}{supp}
\DeclareMathOperator{\maxsupp}{maxsupp}
\DeclareMathOperator{\minsupp}{minsupp}
\DeclareMathOperator{\suc}{succ}

\DeclareMathOperator{\ord}{ord}
\DeclareMathOperator{\cha}{char}

\begin{document}

\title{Strictly singular operators in asymptotic $\ell_p$ Banach spaces}

\author{Anna  Pelczar-Barwacz}

\address{Institute of Mathematics, Jagiellonian University, {\L}ojasiewicza 6, 30-348 Krak\'ow, Poland}

\email{anna.pelczar@im.uj.edu.pl}

\thanks{The research supported by the Polish Ministry of Science and Higher Education grant N N201 421739}


\begin{abstract}
We present condition on higher order asymptotic behaviour of basic sequences in a Banach space ensuring the existence of bounded non-compact strictly singular operator on a subspace. We apply it in asymptotic $\ell_p$ spaces, $1\leq p<\infty$, in particular in convexified mixed Tsi\-rel\-son spaces and related asymptotic $\ell_p$ HI spaces. 
\end{abstract}

\maketitle

\section*{Introduction}
The research on conditions ensuring the existence of non-trivial strictly singular operators on/in Banach spaces increased in last years, in connection with the famous "scalar-plus-compact" problem and following constructions of spaces with "few operators". The "scalar-plus-compact" problem asks if there is an infinite dimensional Banach space on which any bounded operator is a compact perturbation of a multiple of identity. An important step towards solving this problem was made by W.T.~Gowers and B.~Maurey \cite{gm}, who constructed the first HI (hereditarily indecomposable) space, $X_{GM}$, i.e. a space without closed infinite dimensional subspaces which can be written as a direct sum of its further closed infinite dimensional subspaces. Moreover, any operator on a subspace of $X_{GM}$ is a strictly singular perturbation of an inclusion operator. An operator between Banach spaces is \textit{strictly singular}, if none of its restrictions to an infinite dimensional subspace is an isomorphism. The construction of $X_{GM}$ was followed by a class of asymptotic $\ell_1$ HI spaces, started with $X_{AD}$ by S.A.~Argyros and I.~Deliyanni \cite{ad2}, and by a class of asymptotic $\ell_p$ HI spaces \cite{ab,dm}. However, $X_{GM}$ was shown to admit bounded strictly singular non-compact operators first on a subspace \cite{g}, and later - on the whole space \cite{as2}. Also \cite{gas2,b} gave some conditions on parameters of the constructed asymptotic $\ell_p$ HI spaces, ensuring the existence of non-trivial strictly singular operators on the space. Finally the "scalar-plus-compact" problem was solved positively by S.A.~Argyros and R.~Haydon \cite{ah} in the celebrated construction of an HI $\mathscr{L}_\infty$-space with "very few operators". 

A hereditary version of the "scalar-plus-compact" problem, concerning operators on infinite dimensional subspaces of a given space, remains open. Construction of non-trivial strictly singular operators in a Banach space $X$ is based usually on different types of asymptotic behaviour of basic sequences in $X$ with respect to an auxiliary basic sequence $(e_n)$: local representation of $(e_n)$ in $X$, provided for example by Krivine theorem in Lemberg's version \cite{l}, on one side, and "strong" domination of a spreading model of some basic sequence in $X$ by $(e_n)$ on the other \cite{aost,s2,asan}, which ensures strict singularity of the constructed operator. In case of $(e_n)$ equal to the usual basis of $\ell_1$ the asymptotic "strong" domination appears whenever $X$ contains a weakly null basic sequence not generating $\ell_1$-spreading model  \cite{aost}. Construction of non-trivial strictly singular operators based on the higher order representability of $\ell_1$ in a space was studied in \cite{s2}. The operators on the whole space demands specific asymptotic properties of basic sequences in the dual space \cite{as2,gas2,b}. In the last two cases strict singularity is related closely to the hereditary indecomposability of the considered space. 

We present in this paper a general criterium (Th. \ref{main}) ensuring the existence of non-trivial operators in a Banach space in terms of higher order asymptotic behaviour of basic sequences with respect to an auxiliary basic sequence with some regularity properties, under partial unconditionality assumptions. To this end we introduce and study $\al$-strong domination, extending to higher order Schreier families  the notion used in \cite{s2,asan}. Next we apply the general construction in case of any asymptotic $\ell_p$ space $X$ (Cor. \ref{main-as}), providing, as a counterpart of Krivine theorem, "local" lower estimates of basic sequences in $X$ by the usual basis of the $p$-convexified Tsi\-rel\-son-type space $T^{(p)}[\ms_1,\theta]$ with $\theta$ related to asymptotic constants of  $X$ (Th. \ref{theta3}). The further application brings non-trivial strictly singular operators on subspaces of convexified mixed Tsi\-rel\-son spaces and asymptotic $\ell_p$ HI spaces of types constructed in \cite{ab,dm} under with mild conditions on parameters defining the spaces (Cor. \ref{main-as}, \ref{main-xd}). 

The paper is organized as follows: in Section 1. we recall basic notions, in Section 2. we focus on properties of asymptotic $\ell_p$ spaces, proving the "local" lower Tsi\-rel\-son-type estimates. Section 3. is devoted to the study of $\al$-strong domination, for limit $\al<\omega_1$, and in Section 4. we apply developed tools to construct non-trivial operators in general setting and in asymptotic $\ell_p$ spaces, with application to convexified mixed Tsi\-rel\-son spaces and HI spaces.

\section{Preliminaries}
We recall  the basic definitions and standard notation. By a {\em  tree}  we shall mean a non-empty partially ordered  set $(\mt, \preceq)$ such that any set of the form $\{ y \in \mt:y \preceq x \}$, $x \in \mt$, is linearly ordered and finite. If $\mt' \subseteq \mt$ then we say that $(\mt', \preceq)$ is a {\em subtree}  of $(\mt,\preceq)$. The smallest element of a tree (if it exists) is called its \textit{root}, the maximal elements are called {\em terminal}  nodes of a tree. A {\em branch}  in a tree $\mt$ is a maximal linearly ordered set in $\mt$. The \textit{height} of a finite tree is the maximal length of its branches. The {\em immediate successors}  of $t \in \mt$, denoted by $\suc (t)$, are all  the nodes $s \in \mt$ such that $t \preceq s$ but there is no $r \in \mt$ with $t \preceq r \preceq s$. An \textit{order} of a node $t$ of the tree with a root is defined as $\ord (t)=\#\{ s\in\mt: \ s\preceq t\}$.

For any $J\subset\N$ by $[J]^{<\infty}$ we denote the family of finite subsets of $J$. A family $\mf\subset [\N]^{<\infty}$ is \textit{regular}, if it is \textit{hereditary}, i.e. for any $G\subset F$, $F\in \mf$ also $G\in \mf$, \textit{spreading}, i.e. for any integers $n_1<\dots<n_k$ and $m_1<\dots<m_k$ with $n_i\leq m_i$, $i=1,\dots, k$, if $(n_1,\dots,n_k)\in \mf$ then also $(m_1,\dots,m_k)\in \mf$, and \textit{compact} in the product topology of $2^\N$. 

Let $\mf$ be a countable compact family of finite subset of $\N$ endowed with the product topology of $2^\N$.  For any ordinal $\al$ we set $\mf^{\al+1}=\{F\in \mf: F \mbox{\,- a limit point of } \mf^{\al}\}$ and for any limit ordinal $\al$ we set $\mf^{\al}=\cap_{\beta<\al}\mf^{\beta}$. The Cantor-Bendixson index of $\mf$, denoted by $CB(\mf)$, is defined as the least $\al$ for which $\mf^{\al}=\emptyset$. 

\textit{Schreier families} $(\ms_\al)_{\al<\omega_1}$, introduced in \cite{aa}, are defined by induction:
\begin{align*}
\ms_0 &=\{\{ k\}:\ k\in\N\}\cup\{\emptyset\}, \\
\ms_{\al+1}&  =\{F_1\cup\dots\cup F_k:\ k\leq F_1<\dots<F_k, \
F_1,\dots, F_k\in \ms_\al\}, \ \ \al<\omega_1\,.
\end{align*}
If $\al$ is a limit ordinal, choose $\al_n\nearrow \al$ and set
$$
\ms_\al=\{F:\ F\in \ms_{\al_n}\ \mathrm{and}\ n\leq F\ \mathrm{for\ some}\ n\in\N\}\,.
$$
It is well known that the Schreier families $\mc{S}_\al$, $\al<\omega_{1}$, are regular and $CB(\mc{S}_\al)=\omega^\al+1$, $\al<\omega_1$ (c.f \cite{aa}). For any regular family $\mf$ let 
$$
\ms_1(\mf)=\{F_1\cup\dots\cup F_k:\ k\leq  F_1,\dots, F_k\in \mf, \ F_1,\dots, F_k \text{ pairwise disjoint}\}\,.
$$
By an easy adaptation of argument in Lemma 2.1 \cite{m} one can show that $\ms_1(\ms_\al)=\ms_{\al+1}$, $\al<\omega_1$ (cf. also \cite{adkm}). We write $E<F$, for $E,F\subset\N$, if $\max E<\min F$. We say that a sequence $E_1,\dots, E_k$ of subsets of $\N$ is $\ms_\al$-\textit{ad\-missi\-ble}, $\al<\omega_1$, if $E_1<\dots<E_k$ and $(\min E_i)_{i=1}^k\in\ms_\al$.
\bd[$\ms_1$-ad\-missi\-ble tree] 
The $\ms_1$-ad\-missi\-ble tree of finite  subsets of $\N$ is any collection $(E_t)_{t\in\mt}$, indexed by a finite tree $\mt$ with a root 0, such that for any non-terminal node $t\in\mt$ the sequence $(E_s)_{s\in \suc (t)}$ is $\ms_1$-ad\-missi\-ble and  $E_t=\bigcup_{s\in\suc (t)}E_s$.
\ed
\br\label{ad-tree} 
Any $\ms_1$-ad\-missi\-ble tree is a tree ordered by inclusion. By definition of families $(\ms_n)$ for any $\ms_M$-ad\-missi\-ble sequence $(E_k)_k$ of finite subsets of $\N$, $M\in\N$, there is an $\ms_1$-ad\-missi\-ble tree $(E_t)_{t\in\mt}$ of height at most $M$ with $E_0=\bigcup_kE_k$ and $(E_t)_{t\in \mt, t\text{ terminal}}=(E_k)_k$.
\er
Given a Banach space $X$ by $B_X$ denote the closed unit ball of $X$. Let now $X$ be a Banach space with a basis $(e_i)$. The \textit{support} of a vector $x=\sum_i x_i e_i$ is the set $\supp x =\{ i\in \N : x_i\neq 0\}$. We write $x<y$ for vectors $x,y\in X$, if $\supp x<\supp y$. Any sequence $(x_n)\subset X$ with $x_{1}<x_{2}<\dots$ is called a  \textit{block sequence}, a closed subspace spanned by an infinite block sequence $(x_n)$ is called a \textit{block subspace} and denoted by  $[x_n]$. A basic sequence $(x_n)$ $C$-\textit{dominates} a basic sequence $(y_n)$, $C\geq 1$, if for any $(a_n)\in c_{00}$ we have
$$
\nrm{\sum_na_ny_n}\leq C\nrm{\sum_na_nx_n}\,.
$$
Two basic sequences $(x_n)$ and $(y_n)$ are $C$-\textit{equivalent}, $C\geq 1$, if $(x_n)$ $C$-dominates $(y_n)$ and $(y_n)$ $C$-dominates $(x_n)$.
We shall use also the following notion of partial unconditionality \cite{dosz} and equivalence of basic sequences.
\bd 
Let $\mf$ be a family of finite subsets of $\N$.

\cite{dosz} A basic sequence $(x_i)$ is $\mf$-unconditional, if $\norm[\sum_{i\in F}a_ie_i]\leq C\norm[\sum_ia_ie_i]$ for any $(a_i)\in c_{00}$, any $F\in\mf$ and some universal $C\geq 1$. 

We say that basic sequences $(x_i)$ and $(y_i)$ are $\mf$-equivalent, if $(x_i)_{i\in F}$ and $(y_i)_{i\in F}$ are $C$-equivalent for any $F\in\mf$ and some universal $C\geq 1$.
\ed
In the language above a basic sequence $(x_i)$ generates a spreading model $(e_i)$ \cite{bs}, iff for any $\e>0$ for some $n\in\N$ sequences $(e_i)_{i>n}$ and $(x_i)_{i>n}$ are $\ms_1$-equivalent with constant $1+\e$. A basic sequence $(x_i)$ generates an $\ell_1^\al$-spreading model, $\al<\omega_1$ \cite{adm}, iff it is $\ms_\al$-equivalent to the u.v.b. of $\ell_1$.   

We recall that a Banach space $X$ with a basis is \textit{$\ell_p$-asymptotic}, $1\leq p\leq\infty$, if any normalized block sequence $n\leq x_1<\dots<x_n$ is $C$-equivalent to the u.v.b. of $\ell_p^n$, for any $n\in\N$ and some universal $C\geq 1$. 

Finally we say that a sequence $x_1<\dots <x_n$ is $\ms_\al$-\textit{ad\-missi\-ble}, $\al<\omega_1$, if $(\supp x_i)_{i=1}^n$ is $\ms_\al$-ad\-missi\-ble. 
\bd[$p$-convexified mixed Tsi\-rel\-son space]\cite{dm} 
Fix $1\leq p<\infty$, a set $N\subset\N$ and scalars $(\theta_n)_{n\in N}\subset (0,1)$. Define a norm $\nrm{\cdot}$ on $c_{00}$ as the unique norm on $c_{00}$ satisfying the equation
$$
\nrm{x}=\max\left\{\nrm{x}_\infty,\sup\left\{\theta_n^{1/p}(\sum_i\nrm{E_ix}^p)^{1/p}: (E_i)\ \ms_{n}\text{-ad\-missi\-ble},  n\in N\right\}\right\}
$$
The \textit{$p$-convexified mixed Tsi\-rel\-son space} $T^{(p)}[(\ms_n,\theta_n)_{n\in N}]$ is the completion of $(c_{00}, \nrm{\cdot})$. 
\ed 
Take $1<q\leq \infty$ with $\frac{1}{p}+\frac{1}{q}=1$. It is standard to verify that $\nrm{x}=\sup\{f(x):f\in K\}$, $x\in c_{00}$, where $K\subset c_{00}$ is the smallest set such that
\bnum
 \item[(K1)] $(\pm e_i^*)_i\subset K$,
 \item[(K2)] for any  $\ms_{n}$-ad\-missi\-ble $(f_i)\subset K$, $n\in N$, and any $(\gamma_i)\in B_{\ell_q}$ we have $\theta_n^{1/p}\sum_i\gamma_if_i\in K$.
\enum 
In case $p=1$ we obtain the classical mixed Tsi\-rel\-son space $T[(\ms_{n}, \theta_n)_{n\in N}]$, introduced in \cite{ad2}. Notice that for any $p>1$ the space $T^{(p)}[(\ms_{n},\theta_n)_{n\in N}]$ is the $p$-convexification of $T[(\ms_{n}, \theta_n)_{n\in N}]$ \cite{dm} and is  $\ell_p$-asymptotic. It follows immediately by the definition of the space that the u.v.b. $(e_n)$ is 1-unconditional in $T^{(p)}[(\ms_{n},\theta_n)_{n\in N}]$.  

If $N=\{n\}$, we obtain the classical $p$-convexified Tsi\-rel\-son-type space $T^{(p)}[\ms_n,\theta]$. The space $T[\ms_1,1/2]$ is the famous Tsi\-rel\-son space. For $\theta=1$ we have $T^{(p)}[\ms_n,1]=\ell_p$.
We will shorten the notation by denoting any space $T^{(p)}[\ms_1,\theta]$ by $T^{(p)}_\theta$. We recall Lemma 4.13 \cite{otw}: for any sequence $(\theta_n)\subset (0,1]$, with  $\theta_{n+m}\geq \theta_n\theta_m$, $n,m\in\N$,  $\lim_{n\to\infty}\theta_n^{1/n}$ exists and is equal to $\sup_n\theta_n^{1/n}$.
\bn 
A space $T^{(p)}[(\ms_n,\theta_n)_{n\in\N}]$ with $\theta_n\searrow 0$ and $\theta_{n+m}\geq\theta_n\theta_m$  is called a \textit{regular} space. In this case we define  $\theta=\lim_n\theta_n^{1/n}\in (0,1]$.
\en
\br\label{reg}
It follows straightforward that any convexified mixed Tsi\-rel\-son space $T^{(p)}[(\ms_{n},\theta_n)_{n\in N}]$, with infinite $N\subset \N$ and $\theta_n\to 0$, is isometric to a regular space  $T^{(p)}[(\ms_n,\bar{\theta}_n)_{n\in\N}]$, with $\bar{\theta}_n=\sup\{\prod_{i=1}^l\theta_{n_i}:\ \sum_{i=1}^ln_i\geq n, \ n_1,\dots, n_l\in N\}$, $n\in\N$. 
\er
The following notion provides a useful tool for estimating norms in convexified mixed Tsi\-rel\-son spaces.
\bd
[The tree-analysis of a norming functional]\label{def-tree} 
Let $f\in K$, where $K$ is the norming set of a convexified mixed Tsi\-rel\-son space $T^{(p)}[(\ms_n,\theta_n)_{n\in N}]$. By a \textit{tree-analysis} of $f$ we mean a finite family $(f_t)_{t\in \mt}$ indexed by a tree $\mt$ with a unique root $0\in \mt$ satisfying the following:
\bnum
 \item $f_0=f$ and $f_t\in K$ for all $t\in \mt$,
 \item $t\in \mt$ is terminal if and only if $f_t\in (\pm e_n^*)$,
 \item for any non-terminal $t\in \mt$ there is some $n\in N$ such that $(f_s)_{s\in\suc (t)}$ is an $\ms_n$-ad\-missi\-ble sequence and $f_t=\theta_n^{1/p}(\sum_{s\in\suc (t)}\gamma_sf_s)$ for some $(\gamma_s)_{s\in \suc(t)}\in B_{\ell_q}\setminus \{0\}$. In such a case the character of $f_t$ is defined as $\cha(f_t)=n$. 
\enum
\ed
Notice that any $f\in K$ admits a tree-analysis, not necessarily unique.

\section{Lower Tsi\-rel\-son-type estimate in asymptotic $\ell_p$ spaces}
Throughout this section we assume that  $X$ is an asymptotic $\ell_p$ space, $1\leq p\leq \infty$, with a basis. 

For any $n\in\N$ define the lower asymptotic constant $\theta_n=\theta_n(X)\in (0,1]$ (in case $p=1$ cf. \cite{otw})  as the biggest constant such that for any $\ms_n$-ad\-missi\-ble block sequence $n\leq x_1<\dots < x_k\in X$ we have $\norm[x_1+\dots+x_k]^p\geq \theta_n(\norm[x_1]^p+\dots+\norm[x_k]^p)$. It follows easily that $\theta_{n+m}\geq \theta_n\theta_m$, $n,m\in\N$.  Let $\theta=\lim_n\theta_n^{1/n}\in (0,1]$. We will not make at this point the standard stabilization of the constants over block subspaces, or tail subspaces, as it will be done later to satisfy more restrictive conditions. 

The model space for the above situation is a  regular convexified mixed Tsi\-rel\-son space $T^{(p)}[(\ms_n,\theta_n)_{n\in\N}]$. Indeed, by the Fact \ref{est-basis} below  and the definition of the space $(\theta_n)$ is the sequence of its lower asymptotic constants. 
\bfa\label{est-basis} 
Let $Z=T^{(p)}[(\ms_n,\theta_n)_{n\in\N}]$ be a regular $p$-convexified mixed Tsi\-rel\-son space. Then for any $n\in\N$ and $\de>0$ there is a vector $x=\sum_{i\in F}a_ie_i$ with $F\in\ms_n$ such that $\norm[x]\leq(\theta_n^{1/p}+\e)(\sum_{i\in F}|a_i|^p)^{1/p}$. 
\efa
\bp 
By Lemma 1.6 \cite{ad2} for any $n\in\N$ and $\de>0$ there is $(b_i)_{i\in F}\subset (0,1)$, $F\in\ms_n$, such that $\sum_{i\in F}b_i=1$ and $\sum_{i\in G}b_i<\de$ for any $G\in\ms_{n-1}$.
Let $x=\sum_{i\in F}b_i^{1/p}e_i$. Take a norming functional $f\in K$ with a tree-analysis $(f_t)_{t\in\mt}$ and let $G$ be the set of all terminal nodes of $\mt$ with order smaller than $n$. Then $G\in\ms_{n-1}$ and by H\"older inequality and regularity of $(\theta_n)$
\begin{align*}
f(x)&=f(\sum_{i\in G\cap F}b_i^{1/p}e_i)+f(\sum_{i\in F\setminus G}b_i^{1/p}e_i)\\
&\leq (\sum_{i\in G\cap F}b_i)^{1/p}+\theta_n^{1/p}(\sum_{F\setminus G}b_i)^{1/p}<\de^{1/p}+\theta_n^{1/p}\,.
\end{align*}
\ep
In the sequel we will generalize some of the estimates known for $Z$  \cite{kmp} to the case of arbitrary asymptotic $\ell_p$ space $X$. The following Theorem generalizes Lemma 2.14 \cite{kmp} (in case of mixed Tsi\-rel\-son spaces) and Prop. 3.3 \cite{adm} (in case of $\theta=1$), providing also block sequences with supports of uniformly bounded admissibility. One can view this result in context of Krivine theorem in Lemberg's version \cite{l}, stating that for any basic sequence $(x_i)$ there is some $1\leq p\leq \infty$, such that for any $M\in\N$ and $\delta>0$ there is a block sequence $(x_i^{(n)})$ such that any its subsequence of length $M$ is $(1+\de)$-equivalent to the u.v.b. of $\ell_p$. In case of asymptotic $\ell_p$ spaces we increase the order of sequences uniformly "representing" (more precisely dominating) the u.v.b. of some $T^{(p)}[\ms_1,\theta]$ from sequences of fixed length to $\ms_M$-ad\-missi\-ble.

\bt\label{theta3}
Let $X$ be an asymptotic $\ell_p$ space, $1\leq p<\infty$, with lower asymptotic constants $(\theta_n)$. Let $\theta=\lim_n\theta_n^{1/n}$. Then for every  $M\in\N$ and $\delta>0$, there is a normalized block sequence  $(x_i)\subset X$ satisfying for any $G\in \ms_{M}$ and scalars $(a_{i})_{i\in G}$
\begin{equation*}\label{thetae}
\norm[\sum_{i\in G}a_{i}x_{i}] \geq \frac{1}{2}(1-\delta) \norm[\sum_{i\in G}a_{i}e_{\minsupp x_i}]_{T^{(p)}_\theta}.
\end{equation*}
Moreover $(x_i)$ can be chosen to satisfy $(\supp x_i)_i\subset \ms_r$ for  some $r\in\N$. 
\et
In order to achieve the "Moreover..." statement in the above Proposition we introduce more precise lower asymptotic constants measuring the asymptoticity on block sequences with supports of the same admissibility. 

For any normalized block sequence $\textbf{x}=(x_i)\subset X$ and any $n\in\N$ let $\widetilde{\eta}_n(\textbf{x})\in (0,1]$ be the biggest constant such that for any $\ms_n$-ad\-missi\-ble block subsequence $x_{i_1}<\dots < x_{i_k}$ and any scalars $(a_i)_{i=1}^k$ we have $\norm[a_1x_{i_1}+\dots+a_kx_{i_k}]^p\geq \widetilde{\eta}_n(\textbf{x})(|a_1|^p+\dots+|a_i|^p)$. Then let 
$$
\eta_n(\textbf{x})=\sup_{k\in\N}\widetilde{\eta}_n((x_i)_{i\geq k})
$$ 
and finally for any $n\in\N$ let
\begin{align*}
\eta_n=\inf\{\eta_n(\textbf{x}):&\ \textbf{x}=(x_i) \text{ - a normalized block sequence}\\ 
&\text{with }(\supp x_i)\subset\ms_{r_\textbf{x}} \text{ for some }r_{\textbf{x}}\in\N \}\,. 
\end{align*}
It is clear  that $\eta_{n+m}\geq \eta_n\eta_m$, $n,m\in\N$.  Let $\eta=\lim_n\eta_n^{1/n}\in (0,1]$. As $\eta_n\geq \theta_n$ for any $n\in\N$ we have also $\eta\geq \theta$, therefore it will be sufficient to prove the estimate in Th. \ref{theta3} for $T_\eta^{(p)}$ instead of $T_\theta^{(p)}$. 

The proof of  Th. \ref{theta3} is based on the following facts.
\bl\label{sp-vec}
For any $M\in\N$ there is a block sequence $(x_i)\subset X$ such that for any $1\leq j< M$ there is some $\ms_j$-ad\-missi\-ble $(E_k)$ with $\norm[x_i]^p\leq 2\eta^j\sum_k\norm[E_kx_i]^p$, $i\in \N$, and $(\supp x_i)\subset\ms_r$ for some $r\in\N$.
\el
\bp Notice first that for any $M\in\N$ we have
$$
(\sqrt[m]{\eta_m})^M\leq \sqrt[m]{\eta_{Mm}}\leq\sqrt[m]{\eta^{mM}}\,,
$$
thus $\lim_{m\to\infty}\sqrt[m]{\eta_{Mm}}=\eta^M$. Fix $M\in\N$ and by the above pick  $m\in\N$ such that $2^{1/m}\eta_{mM}^{1/m}< 2\theta_M^{1/m}\eta^M$. By definition of $\eta_{mM}$ pick a block sequence $(y_i)\subset X$ with $\norm[y_i]^p\leq 2\eta_{mM}\sum_{k}\norm[F_ky_i]^p$ for some $\ms_{mM}$-ad\-missi\-ble $(F_k)$ and $(\supp y_i)\in\ms_r$ for some $r\in\N$. 

Fix $i\in\N$, let $y=y_i$ and assume that for any $z\in X$ with $\supp z\subset \supp y$ there is some $1\leq j<M$ such that  $\norm[z]^p> 2\eta^j\sum_k\norm[E_kz]^p$ for any $\ms_j$-ad\-missi\-ble $(E_k)$. Notice that if we arrive to contradiction, as $i\in\N$ is arbitrary, we will finish the proof of Lemma. 

Take an $\ms_1$-ad\-missi\-ble tree  $(F_t)_{t\in\mt}$ associated to $(F_k)_k$ as in Remark \ref{ad-tree}.
We will choose inductively some subtree $\mc{R}\subset\mt$ with the same root such that
\bnum
\item $\ord_\mt (t)> (m-1)M$ for any terminal $t\in\mc{R}$,
\item if $t\in\mc{R}$ is non-terminal, then for some $1\leq j_t\leq M$ the sequence $(F_s)_{s\in\suc_{\mc{R}}(t)}$ is $\ms_{j_t}$-ad\-missi\-ble and  $\norm[F_ty]^p\geq 2\eta^{j_t}\sum_{s\in\suc_{\mc{R}}(t)}\norm[F_sy]^p$.
\enum
Notice first that length of the branch linking any terminal node $t$ of $\mc{R}$ and the root is at least $m$  and $\norm[F_ty]^p\geq \theta_M\sum_{F_k\subset F_t}\norm[F_ky]^p$ as  $(F_k)_{F_k\subset F_t}$ is $\ms_{mM-\ord_\mt(t)}$-ad\-missi\-ble and thus also $\ms_M$-ad\-missi\-ble by (1).
Therefore
\begin{align*}
 2\eta_{mM}\sum_{k}\norm[F_ky]^p\geq \norm[y]^p&\geq 2^{m}\sum_{t\in\mc{R}, t\text{ terminal}}\eta^{\ord_\mt(t)}\norm[F_ty]^p\\
&\geq 2^{m}\sum_{t\in\mc{R}, t\text{ terminal}}\eta^{\ord_\mt(t)} \theta_M\sum_{F_k\subset F_t}\norm[F_ky]^p\\
&\geq 2^{m}\eta^{mM}\theta_M\sum_k\norm[F_ky]^p\,,
\end{align*}
hence $2\eta_{mM}\geq 2^m\theta_M\eta^{mM}$ which contradicts the choice of $m$.

We proceed to definition of the tree $\mc{R}$. By our assumption on $y$, considering $z=y$ we have  $\norm[y]^p\geq 2\eta^{j_0}\sum_{s\in\mt, \ord (s)=j_0}\norm[F_sy]^p$ for some $1\leq j_0 \leq M$. Let $\suc_\mc{R}(0)=\{ s\in\mt, \ord_\mt(s)=j_0\}$. Assume we have defined $t\in\mc{R}$ with order $\leq (m-1)M$. By our assumption on $y$, considering $z=F_ty$ we can pick some $1\leq j_t\leq M$ with $\norm[F_ty]^p\geq 2\eta^{j_t}\sum_{s\in\mt, \ord_\mt(s)=\ord_\mt(t)+j_t, F_s\subset F_t}\norm[F_sy]^p$. Let $\suc_{\mc{R}}(t)=\{ s\in\mt, \ord_\mt(s)=\ord_\mt(t)+j_t, F_s\subset F_t\}$ and thus we finish the construction of $\mc{R}$.
\ep
\bfa\label{fact} For any $G\in\ms_M$ and any $z=\sum_{i\in G}c_ie_i\in T_\eta^{(p)}$ there is an $\ms_1$-ad\-missi\-ble tree $\mr$ of height at most $M$, with terminal nodes $\{i\}$, $i\in F$ for some $F\subset G$, of orders $(l_i)_{i\in F}\subset \{1,\dots, M\}$ satisfying  $\norm[z]_{T_\eta^{(p)}}^p\leq 2^p\sum_{i\in F}\eta^{l_i}|c_i|^p$.
\efa
\bp
Take a norming functional $g=\sum_{i\in G}\eta^{k_i/p}\gamma_ie_i^*$ with $(\gamma_i)_{i\in G}\in B_{\ell_q}$ and tree-analysis $(g_t)_{t\in\mt}$ satisfying $g(z)=\norm[z]_{T_\eta^{(p)}}$. Let $I=\{i\in G:\ k_i\leq M\}$.
Let  $g_1$ be the restriction of $g$ to $I$ and $g_2=g-g_1$. If $g_1(z)\geq g_2(z)$ then
$$
g(z)\leq 2g_1(z)\leq 2\sum_{i\in I}\eta^{r_i/p}|\gamma_ic_i| \leq 2(\sum_{i\in I}\eta^{k_i}|c_i|^p)^{1/p}\,,
$$
and we take the tree $\mr=(\supp g_t\cap I)_{t\in\mc{R}}$.
If $g_1(z)\leq g_2(z)$ compute
$$
g(z)\leq 2g_2(z)\leq 2\eta^{M/p}\sum_{i\in G\setminus I}|\gamma_ic_i|\leq 2\eta^{M/p}(\sum_{i\in G}|c_i|^p)^{1/p}\,,
$$
and we take a tree $\mr$ associated to $\ms_M$-ad\-missi\-ble $(\{i\})_{i\in G}$ by Remark \ref{ad-tree}.
\ep
\bp[Proof of Th. \ref{theta3}]
The proof follows the idea of the  proof of Lemma 2.14 \cite{kmp}. Assume the contrary. As in the proof of Lemma \ref{sp-vec} for any $M\in\N$ we have $\lim_{m\to\infty}\sqrt[m]{\eta_{Mm}}=\eta^M$ . Pick $m\in\N$ such that $\eta_{Mm}^{1/m}>2^{1/m}(1-\delta)^p\eta^M$. Take a block sequence $(x_i^0)_i$  according to Lemma \ref{sp-vec} for $mM\in\N$, with $(\supp x^0_i)\subset\ms_r$, for some  $r\in\N$. 

Since the assertion fails there is an infinite sequence $(G_k^{1})_k$ of successive elements of $\ms_{M}$ and coefficients $(a_i^{1})_{i\in G_k^{1},k}$ such that
$$
\norm[\sum_{i\in G_k^{1}}a_i^{1}x_i^0]<\frac{1}{2}(1-\delta) \norm[\sum_{i\in G_k^1}a^{1}_i\norm[x_i^0]e_{m^0_i}]_{T^{(p)}_\eta},\text { for each }k\in\N\,,
$$
where $m_i^0=\minsupp x_i^0$ for each $i$. For any $k\in\N$ set $x_k^{1}=\sum_{i\in G_k^{1}}a_i^{1}x^0_i$ and by Fact \ref{fact} take an $\ms_1$-ad\-missi\-ble tree $\mr_k^1$ with the root $F_k^1\subset G_k^1$ and terminal nodes $(\{i\})_{i\in F_k^1}$, $F_k^1\subset G_k^1$, of orders $(l_i^1)_{i\in F^1_k}\subset\{1,\dots, M\}$ satisfying
$$
\norm[\sum_{i\in G_k^1 } a^{1}_i\norm[x_i^0]e_{m_i^0}]^p_{T^{(p)}_\eta} \leq 2^p\sum_{i\in F_k^1}\eta^{l_i^1}|a^{1}_i|^p\norm[x_i^0]^p \,.
$$
Assume that we have defined  $(x_k^{j-1})_k$ and $(\mr_k^{j-1})_k$ with terminal nodes of orders $(l_i^{j-1})_{i\in F_k^{j-1},k}$ for some $j\leq m$. Then the failure of the assertion  implies the existence of a sequence $(G_k^j)_k$ of successive elements of $\ms_{M}$ and a sequence $(a_i^j)_{i\in G_k^j,k}$ such that for any $k\in\N$
$$
\norm[\sum_{i\in G_k^j}a_i^jx^{j-1}_i]<\frac{1}{2}(1-\delta) \norm[\sum_{i\in G_k^j}a^j_i \norm[x^{j-1}_i]e_{m_i^{j-1}}]_{T^{(p)}_\eta}\,,
$$
where $m_i^{j-1}=\minsupp x_i^{j-1}$ for each $i$. For any $k\in\N$ set $x_k^j=\sum_{i\in G_k^j}a_i^jx^{j-1}_i$ and by Fact \ref{fact} take an $\ms_1$-ad\-missi\-ble tree $\mr_k^j$ with terminal nodes $(\{i\})_{i\in F_k^j}$, $F_k^j\subset G_k^j$, of orders $(l_i^j)_{i\in F_k^j}\subset\{1,\dots, M\}$ satisfying
$$
\norm[\sum_{i\in G_k^j}a^j_i\norm[x_i^{j-1}]e_{m_i^{j-1}}]^p_{T^{(p)}_\eta}\leq 2^p\sum_{i\in F_k^j} \eta^{l_i^j}|a^j_i|^p\norm[x_i^{j-1}]^p \text { for each }k\in\N\,.
$$
The inductive construction ends once we get sequences  $(x_{k}^{m})_k$ and  $(\mr_k^m)_k$.

By the construction  for any $1\leq j\leq m,k\in\N$ we have
\begin{equation}\label{thetae2}
 \norm[x_k^j]^p<(1-\delta)^p\sum_{i\in G_k^j} \eta^{l^j_i}|a^j_i|^p\norm[x_i^{j-1}]^p\,.
\end{equation}
Put $G_k=\cup_{k_{m-1}\in G_k^m}\cup_{k_{m-2}\in G_{k_{m-1}}^{m-1}}\dots\cup_{k_1\in G_{k_2}^2}G_{k_1}^1$, and analogously define $F_k$, for each $k\in\N$. Fix $k\in\N$ and inductively, beginning from $\mr_{k}^{m}$ produce an $\ms_1$-ad\-missi\-ble tree $\mr_k$ by substituting each terminal node $\{ i\}$ of $\mr_{k_j}^j$, $j=1,\dots,m$, by the tree $\mr_i^{j-1}$. Let $(\{i\})_{i\in F_k}$ be the collection of terminal nodes of $\mr_k$ with orders $(l_i)_{i\in F_k}$. Notice that $l_i\leq mM$ for any $i\in F_k$, as each $l_i^j\leq M$. We compute the norm of $x_k^m$, which is of the form
$$
x_k^m=\sum_{k_{m-1}\in G_k^m}\sum_{k_{m-2}\in G_{k_{m-1}}^{m-1}}\dots\sum_{k_1\in G_{k_2}^2}\sum_{i\in G_{k_1}^1}a_{k_{m-1}}^{m}\dots a_i^1x_i^0=\sum_{i\in G_k}b_ix_i^0\,.
$$ 
By the choice of $(x_i^0)$, for  any $i\in\N$ there is an $\ms_{mM-l_i}$-allowable sequence $(E_l)_{l\in L_i}$ with $\norm[x_i^0]^p\leq 2\eta^{mM-l_i}\sum_{l\in L_i}\norm[E_lx_i^0]^p$.

For each $k\in\N$ we have on one hand by repeated use of \eqref{thetae2} 
\begin{align*}
\norm[x_k^m]^p&\leq (1-\delta)^{pm}\sum_{i\in F_k}\eta^{l_i}b_i^p\norm[x_i^0]^p \\
& \leq (1-\delta)^{pm}2\sum_{i\in F_k}\eta^{l_i}b_i^p\eta^{mM-l_i}\sum_{l\in L_i}\norm[E_lx_i^0]^p\\
&=(1-\delta)^{pm}2\eta^{mM}\sum_{i\in F_k}b_i^p\sum_{l\in L_i}\norm[E_lx_i^0]^p\,.
\end{align*}
On the other hand for each $k\in\N$ the sequence $(E_l)_{l\in L_i,i\in F_k}$ is $\ms_{mM}$-ad\-missi\-ble by the definition of $\mr_k$. Consider the block sequence 
$(E_lx_i^0)_{l\in L_i, i\in F_k, k\in\N}$ and notice that $E_l\cap\supp x_i^0\in\ms_r$, for each $l\in L_i,i\in F_k,k\in\N$, by the choice of $(x_i^0)$. Thus by definition of $\eta_{mM}$ for some $k_0\in\N$ we have  
$$
\norm[x_{k_0}^m]^p\geq \eta_{mM}\sum_{i\in F_{k_0}}b_i^p\sum_{l\in L_i}\norm[E_lx_i^0]^p\,,
$$
which brings $\eta_{mM}\leq (1-\delta)^{pm}2\eta^{mM}$, a contradiction with the choice of $m$.
\ep
\br\label{gen} 
In case of $\ell_p^\al$-asymptotic spaces, $1\leq p<\infty$, $\al<\omega_1$,  where all normalized $\ms_\al$-ad\-missi\-ble sequences are uniformly equivalent to the u.v.b. of $\ell_p$ of suitable size, one can define lower asymptotic constants tested on $\ms_{\al n}$-ad\-missi\-ble sequences (in case $p=1$ studied in \cite{otw}). In this setting one obtains analogous results with Tsi\-rel\-son-type spaces $T^{(p)}[\ms_\al, \theta]$. Since the reasoning in this general case  follows exactly the argument in case $\al=1$ above, just by  replacing families $(\ms_n)$ by $(\ms_{\al n})$, for simplicity we present only this last case.
\er
\section{$\omega$-strong domination}
We examine in this section properties of $\al$-strong domination, a higher order counterpart of "strong domination" in \cite{s2} or "domination on small coefficients" in \cite{asan}.  
Throughout this section we fix a limit ordinal $\al<\omega_1$.

For a pair of  seminormalized basic sequences $(x_i)$,$(y_i)$ consider conditions:

$(\bigstar)$ there  are  countable regular families $(\mf_n)$ on $\N$ with $\mf_n\subset\mf_{n+1}$, $n\in\N$, and $CB(\mf_n)\nearrow\omega^\al$,  such that $\Delta_n\to 0$, where for any $n\in\N$
$$
 \Delta_n=\sup\left\{\norm[\sum_ia_ix_i]:\ \max_{F\in\mf_n}\norm[\sum_{i\in F}a_iy_i]\leq \frac{1}{2^n}, \ \norm[\sum_ia_iy_i]\leq 1, \ (a_i)\in c_{00}\right\}\,.
$$

$(\blacktriangle)$ there are countable regular families $(\mf_n)$ on $\N$ with $\mf_n\subset\mf_{n+1}$, $n\in\N$, and $CB(\mf_n)\nearrow\omega^\al$,  such that for any $(a_i)\in c_{00}$
$$
 \norm[\sum_i a_ix_i]\leq \max_{n\in\N}\frac{1}{2^n}\max_{n\leq F\in\mf_n}\norm[\sum_{i\in F} a_iy_i]\,.
$$
\br \label{mf-ms} Take $(\al_n)$ used to define $\ms_\al$. By Prop. 3.10 \cite{otw} for any $\mf$ with $CB(\mf)< \omega^\al$ there are infinite $J\subset\N$ and $n\in\N$ with $\mf\cap[J]^{<\infty}\subset \ms_n$. Therefore $(\bigstar)$ and $(\blacktriangle)$ imply that for some infinite $J=(j_n)\subset \N$ and $(k_n)\subset\N$, subsequences $(x_i)_{i\in J}$ and $(y_i)_{i\in J}$ satisfy analogous properties with families $(\ms_{k_n}\cap [(j_l)_{l>n}]^{<\infty})$. 
\er
\bd 
Fix two seminormalized basic sequences $(x_i)$,$(y_i)$. We say that $(y_i)$ $\al$-strongly dominates $(x_i)$ if
$(y_i)$ is $\ms_\al$-unconditional, $[y_i]$ does not contain $c_0$ and the pair $(x_i)$, $(y_i)$ satisfies $(\bigstar)$.
\ed
As $\mf_0$ is hereditary and spreading, it contains $\ms_0\cap \{k,k+1,\dots\}$ for some $k$ and thus $\al$-strong domination implies domination. The next observation provides a suitable setting for the above definition by Remark \ref{mf-ms}.
\bfa\label{diamond} Let $(y_i)$ be a seminormalized $\ms_\al$-unconditional basic sequence with $[y_i]$ not containing $c_0$. Then for any $\beta\leq\al$ and $\e>0$,  every block subspace $W\subset [y_i]$ contains a vector $w=\sum_ia_iy_i$ with $\max_{F\in\ms_\beta}\norm[\sum_{i\in F}\pm a_iy_i]<\e\norm[w]$. 
\efa
\bp
We show the Fact by induction on $\beta\leq\al$, following the idea of Lemma 3.6 \cite{otw}. Assume that $(y_i)$ is $\ms_\al$-unconditional with  constant 1. For $n=0$ the statement is obvious. Assume the statement holds for $\gamma<\beta$ for fixed $\beta\leq \al$. 

If $\beta$ is limit, take $(\beta_n)$ used to define $\ms_\beta$ and pick a normalized block sequence  $(z_k)\subset W$, $z_k=\sum_{i\in I_k}a_iy_i$, $k\in \N$,  such that 
$$
\max_{G\in\ms_{\beta_n},G\subset I_k}\norm[\sum_{i\in G}\pm a_iy_i]\leq \frac{1}{2^k}, \ n\leq \max I_{k-1}, \ k\in\N\,.
$$ 
Pick any $F\in\ms_\beta$, then $n\leq F\in\ms_{\beta_n}$ for some $n$. Let $k_0=\min\{k\in\N: \ I_k\cap F\neq\emptyset\}$ and compute, using  $n\leq \maxsupp z_{k_0}$ and the $\ms_\al$-unconditionality (provided $\min I_1$ is big enough to ensure $F\cap I_{k_0}\in\ms_\al$), 
\begin{align*}
 \norm[\sum_{i\in F}\pm a_iy_i]&\leq \norm[\sum_{i\in F\cap I_{k_0}}\pm a_iy_i]+\sum_{k>k_0}\norm[\sum_{i\in F\cap I_k}\pm a_iy_i]\leq 1+\sum_{k>k_0}\frac{1}{2^k}\leq 2\,.
\end{align*}
Consider the family  $A=\{\sum_{k\in L}\pm z_k:\ L\in[\N]^{<\infty}\}$. As $[y_i]$ does not contain $c_0$, $\sup_{w\in A}\norm[w]=\infty$ and thus some $w\in A$ satisfies the desired estimate. 

If $\beta=\gamma+1$, pick a normalized block sequence $(z_k)\subset W$, $z_k=\sum_{i\in I_k}a_iy_i$, $k\in \N$, such that 
$$
\max_{G\in\ms_\gamma,G\subset I_k}\norm[\sum_{i\in G}\pm a_iy_i]\leq 1/(2^k\max I_{k-1}),\ \ \ k\in\N\,.
$$ 
Pick any $F\in\ms_\beta$, write $F$ as  $F=F_1\cup\dots\cup F_m$, for some $m\leq F_1<\dots<F_m\in\ms_\gamma$, let $k_0=\min\{k\in\N: \ I_k\cap F\neq\emptyset\}$ and compute, using the $\ms_\al$-unconditionality (provided $\min I_1$ is big enough to ensure $F\cap I_{k_0}\in\ms_\al$)
\begin{align*}
 \norm[\sum_{i\in F}\pm a_iy_i]&\leq \norm[\sum_{i\in F\cap I_{k_0}}\pm a_iy_i]+\sum_{k>k_0}\sum_{j=1}^m\norm[\sum_{i\in F_j\cap I_k}\pm a_iy_i] \leq 1+\sum_{k>k_0}\frac{1}{2^k}\leq 2\,.
\end{align*}
As in the previous case we obtain a suitable $w\in W$ and finish the proof.
\ep
However the $\al$-strong domination appears to be stronger notion than domination without equivalence, in case of $\ell_1$ the situation is simpler.
\bl 
Let $(x_i)$ be a normalized $\ms_\al$-unconditional basic sequence. Assume no subsequence of $(x_i)$ is $\ms_\al$-equivalent to the u.v.b. of $\ell_1$. Then some subsequence of $(x_i)$ is $\al$-strongly dominated by the u.v.b. of $\ell_1$.
\el
\bp Let $(x_i)$ be $\ms_\al$-unconditional with constant 1. Pick $(\al_n)$ used to define $\ms_\al$.  Assume none of subsequences of $(x_i)$ is $\al$-strongly dominated by the u.v.b. of $\ell_1$. Then there are $\delta>0$ and infinite $L\subset\N$ such that for any infinite $J\subset L$ and any $n\in\N$ there is $k_n>n$ and $(a_i)\in c_{00}(J)$ such that $\max_{F\in\ms_{\al_{k_n}}}\norm[\sum_{i\in F\cap J}a_i]\leq 1/2^{k_n}$, $\sum_i|a_i|\leq 1$ and $\norm[\sum_ia_ix_i]>\de$. 

Let $(x_i^*)$ be the biorthogonal functionals to $(x_i)$. Pick $(a_i)$ as above, by unconditionality assume that $(a_i)\subset (0,1)$. We can assume also that $2^{k_n-2}\de\geq 1$ and $\supp\sum_ia_ix_i>n+1$.  Take $(b_i)\subset [0,1]$ with $\sum_ib_ia_i\geq \de$ and $\norm[\sum_ib_ix_i^*]=1$. Let $G_0=\{i\in J:\ b_i>\frac{\de}{4}\}$. Notice that $G_0\not\in\ms_{\al_{k_n}}$, otherwise we arrive to contradiction by the following
$$
\de\leq \sum_ib_ia_i\leq \sum_{i\not\in G_0}b_ia_i+\sum_{i\in G_0}b_ia_i\leq \frac{\de}{4}+\frac{1}{2^{k_n}}\leq \frac{\de}{2}\,.
$$
Pick any $G_1\subset G_0$ with $G_1\in\ms_{\al_{n+1}}\setminus\ms_{\al_n}$. As $G_0>n+1$, also $G_1>n+1$. For any $(c_i)_{i\in G_1}\subset [0,1]$ we have $\norm[\sum_{i\in G_1}c_ix_i]\geq \sum_{i\in G_1}b_ic_i\geq \frac{\de}{4}\sum_{i\in G_1}c_i$, thus by $\ms_\al$-unconditionality $(x_i)_{i\in G_1}$ is $4/\de$-equivalent to the u.v.b. of $\ell_1^{\# F}$. 

Let $\mg$ be the collection of all finite $G\subset L$ such that $(x_i)_{i\in G}$ is $4/\de$-equivalent to the u.v.b. of $\ell_1^{\# G}$. Obviously $\mg$ is hereditary. By the above $\mg\cap [J]^{<\infty}\not\subset\ms_{\al_n}$ for any infinite $J\subset L$ and any $n\in\N$. Therefore by dichotomy  \cite{gas1} there are $J_0\supset J_1\supset \dots$ with $\ms_{\al_n}\cap [J_n]^{<\infty}\subset \mg$, $n\in\N$. It follows that the subsequence $(x_i)_{i\in N}$, where $N=(\min J_n)$, is $\ms_\al$-equivalent to the u.v.b. of $\ell_1$. 
\ep
A typical example of $\omega$-strong domination is formed by convexified mixed Tsi\-rel\-son spaces and Tsi\-rel\-son-type spaces, as the following observation shows.
\bl\label{tsi-star}
Assume $Z=T^{(p)}[(\ms_n,\theta_n)_{n\in\N}]$ is a regular $p$-convexified mixed Tsi\-rel\-son space with $\theta_n/\theta^n\to 0$, where $\theta=\lim_n\theta_n^{1/n}$.  Then the u.v.b. of $T_\theta^{(p)}$ $\omega$-strongly dominates the u.v.b. of $Z$.
\el
\bp As $T^{(p)}_\theta$ is reflexive, if $\theta<1$ or $p>1$, and $T_1^{(1)}=\ell_1$, in all cases $T^{(p)}_\theta$ does not contain $c_0$. To prove condition $(\bigstar)$ notice first that by definition $\norm[x]_Z\leq\norm[x]_{T^{(p)}_\theta}$ for any $x\in c_{00}$. Pick $(a_i)\in c_{00}$ with $\norm[\sum_ia_ie_i]_{T_\theta^{(p)}}=1$ and $\norm[\sum_{i\in F}a_ie_i]_{T_\theta^{(p)}}\leq \frac{1}{2^n}$ for any $F\in\ms_n$. Let $\norm[\sum_ia_ie_i]_Z=\sum_{i\in I}(\prod_j\theta_{l_{i,j}}^{1/p})\gamma_i|a_i|$ for some $(l_{i,j})\subset\N$ and $(\gamma_i)\in B_{\ell_q}$. Let $l_i=\sum_j l_{i,j}$, for any $i\in I$, and $K=\{i\in I: \ l_i\leq n\}$, notice that $K\in\ms_n$ and  compute, by regularity of $Z$,
\begin{align*}
 \norm[\sum_{i\in I}a_ie_i]_Z&\leq \norm[\sum_{i\in K}a_ie_i]_Z+\sum_{i\in I\setminus K}\theta_{l_i}^{1/p}\gamma_i|a_i| \\
& \leq \frac{1}{2^n}+(\max_{l\geq n}\frac{\theta_l}{\theta^l})^{1/p}\sum_{i\in I\setminus K}\theta^{l_i/p}\gamma_i|a_i|\\
& \leq \frac{1}{2^n}+(\max_{l\geq n}\frac{\theta_l}{\theta^l})^{1/p}\norm[\sum_ia_ie_i]_{T_\theta^{(p)}}\,,
\end{align*}
which by assumption on $(\theta_n)$ shows condition $(\bigstar)$ for $(e_i)$ in  $Z$ and $(e_i)$ in  $T_\theta^{(p)}$ with families $(\ms_n)$.
\ep
Next two Lemmas provide characterization of $\al$-strong domination and its invariance (up to taking subsequences) under $\ms_\al$-equivalence. Their proofs follow the reasoning of Proposition 2.3 and Lemma 2.4 \cite{s2}, however additional technic is needed in order to deal with higher order families. 
\bl\label{char} Fix two seminormalized basic sequences $(x_i)$, $(y_i)$.
Then
 \bnum
 \item if the pair $(x_i)$, $(y_i)$ satisfies $(\blacktriangle)$, then it also satisfies $(\bigstar)$,
 \item if $(y_i)$ is unconditional and $[y_i]$ does not contain uniformly $c_0^n$, and the pair $(x_i)$, $(y_i)$ satisfies $(\bigstar)$, then for some infinite $J\subset\N$ the pair $(x_i)_{i\in J}$, $(y_i)_{i\in J}$ satisfies $(\blacktriangle)$.
\enum
\el
\bp (1) Take $(a_i)\in c_{00}$ with $\norm[\sum_ia_iy_i]=1$ and $\norm[\sum_{i\in F}a_iy_i]\leq \frac{1}{2^{n_0}}$ for any $F\in \mf_{n_0}$ and compute by the condition $(\blacktriangle)$
\begin{align*}
 \norm[\sum_ia_ix_i]& \leq \max_n\frac{1}{2^n} \max_{F\in\mf_n}\norm [\sum_{i\in F}a_iy_i]\\
&\leq \max\left\{\max_{F\in\mf_{n_0}}\norm [\sum_{i\in F}a_iy_i], \frac{1}{2^{n_0}}\max_{n>n_0}\max_{F\in\mf_{n}}\norm[\sum_{i\in F}a_iy_i]\right\}\leq \frac{1}{2^{n_0}}\,.
\end{align*}
(2) We can assume that $(y_i)$  is 1-unconditional and 1-dominates $(x_i)$.  
Pick $(k_n)\subset\N$, $k_n>3(n+2)$, such that $\Delta_{k_n}<1/8^{n+1}$, $n\in\N$, where $(\Delta_n)_n$ satisfies the condition $(\bigstar)$ for $(x_i)$ and $(y_i)$.

Define a seminormalized basic sequence $(w_i)$ by the formula
$$
\norm[\sum_ia_iw_i]=\norm[\sum_ia_ix_i]+\max_n\frac{1}{2^n}\max_{F\in\ms_1(\mf_n)} \norm[\sum_{i\in F}a_iy_i]\,.
$$
It is clear that $(w_i)$ dominates $(x_i)$, $(y_i)$ 2-dominates $(w_i)$ and the pair $(w_i)$, $(y_i)$ satisfies $(\bigstar)$ with $(\ms_1(\mf_n))$ and $(\bar{\Delta}_n)=(\Delta_n+\frac{1}{2^n})$. Hence it is enough to show the implication in (2) for sequences $(w_i)$ and $(y_i)$. 

As $(y_i)$ is unconditional and its span does not contain uniformly $c_0^n$'s, we have $l_n<\infty$ for any $n\in\N$, where $l_n$ is the supremum of all $l\in\N$ such that  for some $(z_1,\dots,z_l)\in [y_i]$ with pairwise disjoint supports we have  $\norm[z_j]> 1/2\cdot 8^n$, $j=1,\dots,l$, and $\norm[z_1+\dots+z_l]\leq 2^n$. It follows by definition of $(w_i)$ that for any $n$ the constant  $4l_n$ dominates the supremum of all $l\in\N$ such that for some vector $w\in [w_i]$ with $\norm[w]=1$ and some pairwise disjoint $(E_1,\dots,E_l)\subset \mf_n$ we have $\norm[E_jw]> 1/8^n$, $j=1,\dots,l$.

Let $j_n=\max\{k_n+1,4l_n\}$, $n\in\N$, and $J=\{j_n:\ n\in\N\}$. Take $(a_i)\in c_{00}(J)$, with $\norm[\sum_ia_iw_i]=1$. We define inductively a partition of $J$ into pairwise disjoint  $(F_n)$  such that  for any $n\in\N$
\bnum
 \item[(F1)] $F_n\cap\{j_n,j_{n+1},\dots\}\in\ms_1(\mf_{k_n})$,
 \item[(F2)]  $\norm[\sum_{i\in G}a_iw_i]\leq 1/8^{k_{n-1}}$ for any $G\subset F_n$ with $G\in\mf_{k_{n-1}}$,
 \item[(F3)] if $F_n\neq\emptyset$, then $F_n$ contains some $F\in\mf_{k_n}$ with $\norm[\sum_{i\in F}a_iw_i]> 1/8^{k_n}$,
 \item[(F4)] $\norm[\sum_{i\in F}a_iw_i]\leq 1/8^{k_n}$ for any $F\cap (F_1\cup\dots\cup F_n)=\emptyset$ with $F\in\mf_{k_n}$.
\enum
The first inductive step is similar to the general step, thus we present only the general case. Assume we have $F_1,\dots, F_{n-1}$ satisfying the above. From $I\setminus (F_1\cup\dots\cup F_{n-1})$ we pick a maximal family of pairwise disjoint sets $(F_n^j)_j\subset\mf_{k_n}$ with $\norm[\sum_{i\in F_n^j}a_iw_i]>1/8^{k_n}$ for each $j$. Let $F_n=\cup_jF_n^j$. It follows that conditions (F3) and (F4) are satisfied. As there can be at most $4l_n\leq j_n$ many $(F_n^j)$'s  we obtain (F1). Finally the condition (F4) for $n-1$ implies (F2) for $n$, which ends the inductive construction. Compute, using (F2)
\begin{align*}	
1=\norm[\sum_ia_iw_i]&\leq\sum_{n}\sum_{i\in F_n,i<j_n}|a_i|+\sum_{n}\norm[\sum_{i\in F_n,i\geq j_n}a_iw_i]\\
&\leq \sum_{n} \frac{n}{8^{k_{n-1}}}+\sum_{n} \norm[\sum_{i\in F_n,i\geq j_n}a_iw_i]\,.
\end{align*}
It follows that $1/2\leq \sum_{n} \norm[\sum_{i\in F_n,i\geq j_n}a_iw_i]$ and thus for some $n_0$ we have
$$
\norm[\sum_{i\in F_{n_0},i\geq j_{n_0}}a_iw_i]\geq \frac{1}{2^{n_0+1}}\,.
$$ 
As $(y_i)$ 2-dominates $(w_i)$ we have $ 1/2^{k_{n_0}}\leq \norm[\sum_{i\in F_{n_0}, i>j_{n_0}}a_iy_i]$. On the other hand by (F2) and definition of $(w_i)$ we have $\norm[\sum_{i\in G}a_iy_i] <1/4^{k_{n_0}}$ for any $G\subset F_{n_0}$ with $G\in\mf_{k_{n_0-1}}$. Therefore by $(\bigstar)$ for $(w_i)$ and $(y_i)$ we obtain
$$
\norm[\sum_{i\in F_{n_0},i\geq j_{n_0}}a_iw_i]\leq \bar{\Delta}_{k_{n_0-1}}\norm[\sum_{i\in F_{n_0},i\geq j_{n_0}}a_iy_i]\,.
$$
Putting the estimates together, by the choice of $(k_n)$ and (F1) we obtain
\begin{align*}
 \norm[\sum_ia_iw_i]=1&\leq 2^{n_0+1}\norm[\sum_{i\in F_{n_0}, i>j_{n_0}}a_iw_i]\leq \frac{1}{2^{n_0}}\norm[\sum_{i\in F_{n_0}, i>j_{n_0}}a_iy_i]\\
 &\leq \frac{1}{2^{n_0}}\max_{n_0\leq F\in\ms_1(\mf_{k_{n_0}})}\norm[\sum_{i\in F\cap J}a_iy_i]\\
\end{align*}
which yields $(\blacktriangle)$ for $(x_i)_{i\in J}$ and $(y_i)_{i\in J}$ with families $(\ms_1(\mf_{k_n})\cap [J]^{<\infty})$.
\ep
\bl\label{equi} Consider seminormalized basic sequences $(x_i)$, $(z_i)$, $(y_i)$ with $(y_i)$ unconditional and $[y_i]$ not containing uniformly $c_0^n$'s. 

Assume $(x_i)$ and $(z_i)$ are $\ms_\al$-equivalent. Then if the pair $(z_i), (y_i)$ satisfies $(\bigstar)$, then for some infinite $J\subset\N$ also $(x_i)_{i\in J}, (y_i)_{i\in J}$ satisfies $(\bigstar)$. 
\el
\bp
We can assume that the basic sequence $(x_i)$ is bimonotone and $(y_i)$ is 1-unconditional. Let $C\geq 1$ be the $\ms_\al$-equivalence of $(x_i), (z_i)$ constant. Take $(\al_n)$ used to define $\ms_\al$. 
Take $(\Delta_n)$  satisfying the condition $(\bigstar)$ for $(z_i), (y_i)$ and pick $(k_n)$, $k_n\geq n$, such that $\sum_{n}\Delta_{k_{n-1}}<\infty$.  By Remark \ref{mf-ms} there is $(t_n)\subset\N$ such that $\mf_{k_n}\cap [(t_i)_{i>n}]^{<\infty}\subset\ms_{\al_{t_n}}$ for each $n\in\N$.

Since $[y_i]$ does not contain uniformly $c_0^n$'s, for any $n$ we have $l_n<\infty$, where $l_n$ is the supremum of all $l\in\N$ such that  for some disjointly supported $z_1,\dots,z_l\in [y_i]$ with $\norm[z_j]> 1/2^{k_n}$, $j=1,\dots,l$, we have  $\norm[z_1+\dots+z_l]\leq 1$.

Pick $J=\{j_n: \ n\in\N\}\subset \{t_n\}$ with $j_n\geq \max\{k_n+1,l_n,t_n+1\}$, $n\in\N$. Take $(a_i)\in c_{00}(J)$, with $\norm[\sum_ia_iy_i]=1$.
As in  the proof of Lemma \ref{char} we define inductively a partition of $J$ into pairwise disjoint  $(F_n)$  such that  for any $n\in\N$
\bnum
 \item[(F1)] $F_n\cap\{j_n,j_{n+1},\dots\}\in\ms_1(\mf_{k_n})\subset \ms_{\al_{t_n}+1}$,
 \item[(F2)]  $\norm[\sum_{i\in G}a_iy_i]\leq 1/2^{k_{n-1}}$ for any $G\subset F_n$ with $G\in\mf_{k_{n-1}}$,
 \item[(F3)] if $F_n\neq\emptyset$, then $F_n$ contains some $F\in\mf_{k_n}$ with $\norm[\sum_{i\in F}a_iy_i]> 1/2^{k_n}$,
 \item[(F4)] $\norm[\sum_{i\in F}a_iy_i]\leq 1/2^{k_n}$ for any $F\cap (F_1\cup\dots\cup F_n)=\emptyset$ with $F\in\mf_{k_n}$.
\enum
Now compute
\begin{align*}
\norm[\sum_ia_ix_i]&\leq\sum_{n: F_n\neq\emptyset}\sum_{i\in F_n,i<j_n}|a_i|+\sum_{n: F_n\neq\emptyset}\norm[\sum_{i\in F_n,i\geq j_n}a_ix_i]\\
&\leq \sum_{n: F_n\neq\emptyset} \frac{n}{2^{k_{n-1}}}+C\sum_{n: F_n\neq\emptyset} \norm[\sum_{i\in F_n,i\geq j_n}a_iz_i]\ \ \ \ \ \text{by (F2) and (F1)}\\
&\leq \sum_{n: F_n\neq\emptyset}\frac{n}{2^{k_{n-1}}}+C\sum_{n: F_n\neq\emptyset} \Delta_{k_{n-1}}\ \ \ \ \ \text{by (F2) and }(\bigstar).
\end{align*}
Fix $n_0\in\N$ and assume additionally that $\norm[\sum_{i\in F}a_iy_i]\leq 1/2^{k_{n_0}}$ for any  $F\in\mf_{k_{n_0}}$. Then by (F3), (F4) and the above computation 
$$
\norm[\sum_ia_ix_i]\leq \sum_{n\geq n_0}\frac{n}{2^{k_{n-1}}}+C\sum_{n\geq n_0}\Delta_{k_{n-1}}\,,
$$ 
thus $(\bigstar)$ for $(x_i)_{i\in J}$, $(y_i)_{i\in J}$ is satisfied with families $(\mf_{k_n}\cap [J]^{<\infty})$.
\ep
\section{Strictly singular non-compact operators}
In this section we apply tools developed in the previous part to give sufficient conditions for existence of non-trivial strictly singular operators. 
We note first a version of Theorem 1.1 \cite{s2} in $\ms_\al$-unconditional setting. 
\bpr\label{s2} Let $(x_i)$ and $(y_i)$ be two seminormalized basic sequences such that $(y_i)$ $\al$-strongly dominates $(x_i)$, for some limit $\al<\omega_1$. 

Then the map $y_i\mapsto x_i$ extends to a bounded non-compact strictly singular operator between $[y_i]$ and $[x_i]$.
\epr
\bp As $(y_i)$ dominates $(x_i)$, the map $y_i\mapsto x_i$ extends to a bounded non-compact operator $T$ between $[x_i]$ and $[y_i]$. To prove the strict singularity use $(\bigstar)$ and Fact \ref{diamond} with Remark \ref{mf-ms}.
\ep
The next theorem will serve as a base for further applications. We build an operator using block sequences with different asymptotic behaviour with respect to an auxiliary basic sequence $(e_i)$. However the situation is analogous to the results in  \cite{aost, asan, kmp}, we work on $(\ms_{\al_n})$-ad\-missi\-ble sequence instead of $(\ma_n)$-ad\-missi\-ble sequences, i.e. sequences of length $n$, $n\in\N$. The sequence $(e_i)$ plays the role of a spreading model in \cite{aost, asan, kmp}, in our setting we require domination of $(e_i)$ by all its subsequences instead of subsymmetry. 
\bt\label{main} 
Let $X$ be a Banach space with an $\ms_\al$-unconditional basis, for limit $\al<\omega_1$. Let $E$ be a Banach space with an unconditional basis $(e_i)$ dominated by all its subsequences, not containing uniformly $c_0^n$'s. Assume that
\bnum
\item $X$ has a normalized basic sequence $(x_i)$ $\al$-strongly dominated by $(e_i)$,
\item for any $\beta<\al$ there exists a normalized block sequence $(x_i^\beta)_i$ with $(\supp x_i^\beta)_i\subset\ms_{r_\beta}$, for some $r_\beta\in\N$, such that $(x_i^\beta)_{i\in F}$ $C$-dominates $(e_i)_{i\in F}$  for any $F\in\ms_\beta$ and universal $C\geq 1$.
\enum
Then  $X$ admits a bounded strictly singular non-compact operator on a subspace.
\et
\br
In case $E=\ell_1$ Theorem above follows by~Th. 1.4, \cite{s2}, as (1) and (2) imply (a) and (b) in Th.~1.4. In case of $E=\ell_1$ partial unconditionality of suitable sequences follows by \cite{amt}.  Comparing to Th. 1.4 \cite{s2} Theorem above can be regarded as an extension of Th. 1.4 in replacing the u.v.b. of $\ell_1$ by other basic sequence, however with the price paid on additional assumptions related to partial unconditionality. Recall that by \cite{o} any normalized weakly null sequence admits an $\ms_1$-unconditional subsequence, and the result was extended in \cite{asan} to special arrays of vectors, but analogous statement does not hold for $\ms_\al$ with $\al>1$.
\er
In the proof the lack of full unconditionality is substituted by $\ms_\al$-un\-con\-di\-tio\-na\-lity and uniform bound on admissibility of supports of each of block sequences $(x_i^{(n)})_i$  in (2). It follows that projections on $[(x_i^\beta)_{i\in F}]$ are bounded uniformly on $F\in\ms_\beta$ provided $\min F$ is big enough and $\beta<\al$. We produce a block sequence $(y_i)$ from sequences $(x_i^\beta)$ in the standard way and show that some subsequences $(x_i)_{i\in J}$ and $(y_i)_{i\in J}$ satisfy $(\bigstar)$ passing through Lemma \ref{char}. Since we cannot assure even $\ms_\al$-unconditionality of $(y_i)$, we need to prove strict singularity of the operator carrying $(y_i)_{i\in J}$ to $(x_i)_{i\in J}$ by hand. 
\bp  
Take $(\al_n)$ used to define $\ms_\al$. We can assume that $X$ does not contain $c_0$ and its basis is $\ms_\al$-uncon\-di\-tio\-nal with constant 1. As $(e_n)$ is dominated by all its subsequences, it is also uniformly dominated by its subsequences, and we assume that the uniform domination constant is 1. By Lemma \ref{char} and Remark \ref{mf-ms} for some infinite $J\subset\N$, $(k_n)\subset\N$, we have, letting   $\mf_n=\ms_{\al_{k_n}}$, 
$$
 \norm[\sum_i a_ix_i]\leq \max_{n\in\N}\frac{1}{4^n}\max_{n\leq F\in\mf_n}\norm[\sum_{i\in F\cap J} a_ie_i], \ (a_i)\in c_{00}(J)\,.
$$
Given $(x_i^{\al_n})_i\subset X$, $n\in\N$, as in (2) let $y_i^{(n)}=x_i^{\al_{k_n}}$ for any $i,n\in\N$. By the assumption on $(e_i)$, passing to subsequences we can assume that $y_1^{(1)}<y_2^{(1)}<y_2^{(2)}<y_3^{(1)}<y_3^{(2)}<y_3^{(3)}<\dots$ and $ r_{\al_{k_n}}+k_n<y_i^{(n)}$ for any $i\geq n$. 
Then
\begin{equation}\label{*}
 \supp\sum_{i\in F}y_i^{(n)}\in\ms_\al \text{ for any }n\leq F\in\mf_n,\  n\in\N\,.
\end{equation}
By choice of $(y_i^{(n)})_i$ we have for any $(a_i)\in c_{00}(J)$ 
 \begin{align*}
 \norm[\sum_ia_ix_i]\leq \max_{n\in\N}\frac{C}{4^n} \max_{n\leq F\in\mf_n}\norm[\sum_{i\in F\cap J}a_i y_i^{(n)}]\,. 
\end{align*}
Let $y_i=\sum_{n=1}^i\frac{1}{2^n}y_i^{(n)}$ for any $i\in I$.  Obviously $(y_i)$ is a seminormalized block sequence. 
Fix now $n_0\in \N$ and continue the above estimation
 \begin{align*}
 \norm[\sum_ia_ix_i]&\leq C\max\left\{\max_{\substack{n\leq n_0\\ n\leq F\in\mf_n}}\norm[\sum_{i\in F\cap J}a_i y_i^{(n)}], \frac{1}{4^{n_0}} \max_{\substack{n>n_0\\ n\leq F\in\mf_n}}\norm[\sum_{i\in F\cap J}a_i y_i^{(n)}]\right\}  \\
&\leq C\max\left\{ \max_{n\leq n_0}\max_{ n\leq F\in\mf_{n_0}}\norm[\sum_{i\in F\cap J}a_i y^{(n)}_i], \frac{1}{2^{n_0}}\norm[\sum_ia_i y_i]\right\}\,,
\end{align*}
where the last inequality follows by \eqref{*} and $\al$-unconditionality of $(x_i)$. Thus the following Claim holds true.
\bcl[\bf{A}] For any $n_0\in\N$ and $(a_i)\in c_{00}(J)$ with $\norm[\sum_ia_iy_i]=1$ we have
$$
 \norm[\sum_{i}a_ix_i]\leq C\max\left\{ \max_{n\leq n_0}\max_{ n\leq F\in\mf_n}\norm[\sum_{i\in F\cap J}a_i y^{(n)}_i], \frac{1}{2^{n_0}}\right\}\,.
$$
\ecl
Taking $n_0=0$ we obtain that $(y_i)_{i\in J}$ dominates $(x_i)_{i\in J}$, thus the mapping $y_i\mapsto x_i$ extends to a bounded non-compact operator $T:[(y_i)_{i\in J}]\to [(x_i)_{i\in J}]$. However we obtain also $(\bigstar)$ for the pair $(x_i)_{i\in J}$, $(y_i)_{i\in J}$, without $\ms_\al$-unconditionality of $(y_i)$ we need to prove strictly singularity of $T$ by hand. First we adapt Fact \ref{diamond} to our setting. 
\bcl[\bf{B}] Given any $n\in\N$ and $\e>0$,  any block subspace $W\subset [y_i]$ contains a further block subspace $V$ such that any $w=\sum_ia_iy_i\in V$ satisfies 
$$
\max_{F\in\mf_n}\norm[\sum_{i\in F}a_iy_i^{(n)}]<\e\norm[w]\,.
$$ 
\ecl
To prove Claim (B) we first show that for any $\e>0$,  $n\in\N$, $\beta\leq\al_{k_n}$, any block subspace $W\subset [y_i]$ contains $w_\e=\sum_ia_iy_i$ satisfying $\max_{F\in\ms_\beta}\norm[\sum_{i\in F}a_iy_i^{(n)}]<\e\norm[w_\e]$. The proof of this statement follows step by step the proof of Fact \ref{diamond}, as we assumed at the beginning that $X$ does not contain $c_0$. We assume that $W\geq n$, estimate $\norm[\sum_{i\in F}\pm a_iy_i^{(n)}]$ instead of $\norm[\sum_{i\in F}\pm a_iy_i]$ and use \eqref{*} to obtain $\norm[\sum_{i\in G}a_iy_i^{(n)}]\leq \norm[\sum_{i\in G}a_iy_i]$ for any $n\leq G\in\mf_n$. Once we have this statement, to complete the proof of Claim (B) let $V=[w_{\e/2^i}]$. 

With  the above two Claims we are ready to prove the strict singularity of $T$. Fix $n_0\in\N$, take any block subspace $W\subset [y_i]$ and using Claim (B) pick inductively block subspaces $W\supset V_{n_0}\supset V_{n_0-1}\supset\dots\supset V_0$ such that for any $w=\sum_ia_iy_i\in V_0$ we have $\max_{F\in\mf_n}\sum_{i\in F}a_iy_i^{(n)}\leq\frac{1}{2^{n_0}}\norm[w]$ for any $n\leq n_0$. Claim (A) ends the proof. 
\ep
The model space $E$ in Th. \ref{main} is the $p$-convexified Tsi\-rel\-son-type space $T_\theta^{(p)}$, for $1\leq p<\infty$ and $\theta\in (0,1]$. As Th. \ref{theta3} yields condition (2) of Th. \ref{main} in case $\al=\omega$ for any asymptotic $\ell_p$ space $X$ with lower asymptotic constants $(\theta_n)$ and $E=T_\theta^{(p)}$, where $\theta=\lim_n\theta_n^{1/n}$, we obtain the following.
\bc\label{main-as} Let $X$ be an asymptotic $\ell_p$ Banach space, $1\leq p<\infty$, with lower asymptotic constants $(\theta_n)$ and an $\ms_\omega$-unconditional basis. 

Assume $X$ contains a normalized basic sequence $(x_i)$ $\omega$-strongly dominated by the u.v.b. of $T_\theta^{(p)}$, where $\theta=\lim_n\theta_n^{1/n}$. 

Then $X$ admits a bounded strictly singular non-compact operator on a subspace.
\ec
By Lemma \ref{tsi-star} the typical space $X$ for the above situation is a regular $p$-convexified mixed Tsi\-rel\-son space $T^{(p)}[(\ms_n,\theta_n)_n]$ with $\theta_n/\theta^n\to 0$, where $\theta=\sup_n\theta_n^{1/n}$.  However, as conditions (1) and (2) of Th. \ref{main} are invariant under $\ms_\omega$-equivalence up to taking subsequences (for (1) use  Lemma  \ref{equi}), a stronger result, requiring only $\ms_\omega$-representation of the regular mixed Tsi\-rel\-son space, holds true. 
\bc\label{main-tsi} Let $X$ be a Banach space with an $\ms_\omega$-unconditional basis $(x_i)$.

Assume the basis $(x_i)$ is $\ms_\omega$-equivalent to the u.v.b. of a  regular $p$-conve\-xi\-fied mixed Tsi\-rel\-son space $T^{(p)}[(\ms_n,\theta_n)_n]$. Assume also that $\theta_n/\theta^n\to 0$, where $\theta=\lim_n\theta_n^{1/n}$.

Then $X$ admits a bounded strictly singular non-compact operator on a subspace.
\ec
\br By Remark \ref{gen} above Corollaries hold for any $\al<\omega_1$, in terms of $\ell_p^\al$-asymptotic spaces, convexified mixed Tsi\-rel\-son spaces $T^{(p)}[(\ms_{\al n},\theta_n)]$ and convexified Tsi\-rel\-son-type spaces $T^{(p)}[\ms_\al,\theta]$. 
\er
We will recall now construction of spaces based on mixed Tsi\-rel\-son spaces, initiated in \cite{ad2}, used for building classes of HI asymptotic $\ell_p$ spaces with different types of properties, see also \cite{adkm,ab,dm}.

Fix $1\leq p<\infty$, let $1<q\leq \infty$ satisfy $1/p+1/q=1$. Fix infinite sets $N,L\subset\N$ (not necessarily disjoint) and scalars $(\theta_n)_{n\in N}, (\rho_l)_{l\in L}\subset (0,1)$ with $\theta_n\to 0$, $\rho_l\to 0$. Assume also the regularity of $(\theta_n)$, i.e. that $\theta_n\geq\prod_{i=1}^l\theta_{n_i}$ for any $n,n_1,\dots,n_l\in N$ with $\sum_{i=1}^ln_i\geq n$. 

Let $\mc{W}=\{(f_1,\dots,f_k):\ f_1<\dots<f_k\in c_{00}(\Q)\cap B_{\ell_q}, k\in\N\} $ and fix an injective function $\sigma:\mc{W}\to N$.
For any $D\subset c_{00}(\Q)$ define
\begin{align*}
D_n=&\{\theta_n^{1/p}\sum_{i=1}^k\gamma_i f_i: \ f_1,\dots,f_k\in D,  \  (f_1,\dots,f_k) \text{ is } \ms_n\text{-ad\-missi\-ble},\\
& \ (\gamma_i)\in B_{\ell_q}\cap c_{00}(\Q), \ k\in\N\},  \ \ \ \ \ \ \ \  n\in N \\
D_l^\sigma=&\{\rho_l^{1/p}\sum_{i=1}^k\gamma_i Ef_i: \ f_1,\dots,f_k\in D, \  (f_1,\dots,f_k) \text{ is } (\sigma,\ms_l)\text{-ad\-missi\-ble},\\
& \ (\gamma_i)\in B_{\ell_q}\cap c_{00}(\Q), \ E\subset\N \text{ interval}, \ k\in\N\},   \ \ \ \ \ \ \ \ \  l\in L
\end{align*}
where a block sequence $(f_1,\dots,f_k)$ is $(\sigma, \ms_l)$-ad\-missi\-ble, if $(f_1,\dots,f_k)$ is $\ms_l$-ad\-missi\-ble, $f_1\in\bigcup_{n\in N}D_n$ and $f_{i+1}\in D_{\sigma(f_1,\dots,f_i)}$ for any $i<k$.

Consider a symmetric set $D\subset c_{00}(\Q)$ such that
\bnum
 \item[(D1)] $(\pm e_n^*)_n\subset D$,
 \item[(D2)] $D\subset\bigcup_{n\in N}D_n\cup\bigcup_{l\in L}D_l^\sigma$,
 \item[(D3)] $D_n\subset D$ for any $n\in N$.
\enum
Define a norm on $c_{00}$ by $\norm[x]_D=\sup\{f(x):f\in D\}$, $x\in c_{00}$, denote by $X_D$ the completion of $(c_{00},\norm[\cdot]_D)$. Obviously the u.v.b. $(e_n)$ is a basis for $X_D$.

It follows that $D\subset K_{N\cup L}$, where $K_{N\cup L}$ is the norming set of the $p$-convexified mixed Tsi\-rel\-son space defined by all pairs $(\ms_n, \theta_n)_{n\in N}\cup (\ms_l,\rho_l)_{l\in L}$, thus each functional in $D$ admits a tree-analysis (Def. \ref{def-tree}). By (D3) also $D\supset K$, where $K$ is the norming set of $T^{(p)}[(\ms_n,\theta_n)_{n\in N}]$.
\bc \label{main-xd}
Let $X_D$ be defined as above. Assume $$
\lim_{n\in N, n\to\infty}\theta_n/\theta^n=0,\ \text{ where }\ \theta=\sup_{n\in N}\theta_n^{1/n}\,.
$$
Then $X_D$ admits a bounded non-compact strictly singular operator on a subspace. 
\ec
\bp It is enough to show that for some $(i_n)_n\subset \N$ the following hold
 \bnum
 \item sequence $(e_{i_n})\subset X_D$ is $\ms_\omega$-unconditional,
 \item  sequences $(e_{i_n})\subset X_D$, $(e_{i_n})\subset T^{(p)}[(\ms_n,\theta_n)_{n\in N}]$ are $\ms_\omega$-equivalent.
 \enum
Indeed, recall that $T^{(p)}[(\ms_n,\theta_n)_{n\in N}]$  is isomorphic to a regular space given by $T^{(p)}[(\ms_n,\bar{\theta}_n)_{n\in\N}]$, with $(\bar{\theta}_n)$ defined as in Remark \ref{reg}. By the regularity of $(\theta_n)_{n\in N}$ we have $\theta_n=\bar{\theta}_n$ for any $n\in N$. Therefore the subspace $[e_{i_n}]$ by Lemma \ref{tsi-star} satisfies the assumption of Cor. \ref{main-tsi}, which ends the proof.

Now we pick  $(i_n)_n\subset\N$ with desired properties. Let $Z=T^{(p)}[(\ms_n,\theta_n)_{n\in N}]$. We denote by $(e_i)$ the u.v.b. both in $X_D$ and $Z$. We will show the following 
\bcl
For any $n\in\N$ there is $i_n\in\N$ such that for any $(a_i)_{i\in F}$ with $F\in\ms_n$ and $F\geq i_n$ we have $\norm[\sum_{i\in F}a_ie_i]_D\leq 4\norm[\sum_{i\in F}a_ie_i]_Z$. 
\ecl
First notice that Claim implies (1) and (2) for $(e_{i_n})$. Indeed, (2) follows straightforward, as $\norm[\sum_ia_ie_i]_D\geq \norm[\sum_ia_ie_i]_Z$ for any $(a_i)\in c_{00}$ by the property $K\subset D$.  Also by Claim for any $(a_i)_{i\in F}$ with $i_n\leq F\in\ms_n$, $n\in\N$, there is a norming functional $f\in Z^*$, therefore also $f\in X_D^*$, with $\supp f\subset F$, such that $\norm[\sum_{i\in F}a_ie_i]_D\leq 4f(\sum_{i\in F}a_ie_i)$ in $X_D$. Thus we obtain (1) for $(e_{i_n})$.

We proceed to proof of Claim. Fix $n\in\N$. Pick $j_n$ such that $\theta_{j}\leq \frac{1}{2^p}\theta_n$ for any $j_n\leq j\in N$ and $\rho_j\leq \frac{1}{2^p}\theta_n$ for any $j_n\leq j\in L$. By injectivity of $\sigma$ there is $i_n$ such that $\sigma(f)>j_n$ for any $f\in \mc{W}$ with $\maxsupp f\geq i_n$. 

Take now any $(a_i)_{i\in F}$, $i_n\leq F\in\ms_n$, with $\norm[\sum_{i\in F}a_ie_i]_D=1$. It follows by (D1) and (D3) that $\theta_n\sum_{i\in F}|a_i|^p\leq 1$. Take a norming functional $f\in D$ with a tree-analysis $(f_t)_{t\in \mt}$ satisfying $f(\sum_{i\in F}a_ie_i)=1$. Let 
$$
I=\{i\in F:\ \cha(f_t)<j_n \text{ for any }t\in\mt \text{ with }f_t(e_i)\neq 0\}\,.
$$ 
Then by H\"older inequality and choice of $j_n$ 
$$
|f(\sum_{i\in F\setminus I}a_ie_i)|\leq \frac{1}{2}\theta_n^{1/p}(\sum_{i\in F}|a_i|^p)^{1/p}\leq \frac{1}{2}\,.
$$ 
Thus $f(\sum_{i\in I}a_ie_i)\geq \frac{1}{2}$. Let $I_1=\{i\in I:\ a_i> 0, \ f(e_i)> 0 \}$ and $I_2=\{i\in I:\ a_i<0, \ f(e_i)<0 \}$. Then either $f(\sum_{i\in I_1}a_ie_i)\geq \frac{1}{4}$ or  $f(\sum_{i\in I_2}a_ie_i)\geq \frac{1}{4}$. Assume the first case holds and let $x=\sum_{i\in I_1}a_ie_i$. Take any $t\in\mt$ with $f_t(x)\neq 0$ and $f_t\in D_l^\sigma$ for some $l\in L$. Then by choice of $i_n$ and $I$ there is at most one $s_t\in\suc (t)$ with  $\supp f_{s_t}\cap I\neq\emptyset$. 

Given any non-terminal $t\in\mt$, with $f_t=\theta_{n_t}^{1/p}\sum_{s\in\suc (t)}\gamma_sf_s$ we let $|f_t|=\theta_{n_t}^{1/p}\sum_{s\in\suc (t)}|\gamma_s|f_s$. Construct a functional $g$ replacing in the tree-analysis $(f_t)_{t\in\mt}$ each $f_t\in D_l^\sigma$ by $|f_{s_t}|$. Then $g\in K$ as every node of the tree-analysis of $g$ belongs to $\bigcup_nD_n$. For  $h$ defined as the restriction of $g$ to $I$ we have $h\in K$ and $h(\sum_{i\in F}a_ie_i)= h(x)\geq f(x)\geq \frac{1}{4}$, which ends the proof of Claim.
\ep
\br Notice that in case $\theta=1$ the sequence $(\bar{\theta}_n)$ defined in Remark \ref{reg} also satisfies $\bar{\theta}=1$, thus the assumption of Cor. \ref{main-tsi} are satisfied. Therefore in this case  we do not  need the regularity of $(\theta_n)_{n\in N}$ in Cor. \ref{main-xd}. 
\er
\bc The HI $\ell_2$-asymptotic Banach space $X_{AB}$  constructed in \cite{ab} and HI asymptotic $\ell_p$ Banach  spaces $X_{(p)}$, $1< p<\infty$, constructed in \cite{dm} admit bounded strictly singular non-compact operator on a subspace.
\ec
\bp To show the Corollary notice that spaces $X_{AB}$  and $X_{(p)}$ are of the form $X_D$ with $N=(n_{2i})$, $L=(n_{2i+1})$, $\theta_{n_{2i}}=\frac{1}{m_{2i}^p}$, $i\in \N$, for suitably chosen $(n_i)$, $(m_i)$, satisfying $\theta_{n_{2i}}^{1/n_{2i}}\to 1$. In case of $X_{AB}$ we have $\rho_{n_{2i+1}}=\frac{1}{m_{2i+1}^2}$, in case of $X_{(p)}$ we have $\rho_{n_{2i+1}}=\frac{2}{2^pm_{2i+1}^p}$. The Remark above ends the proof. 
\ep
\br Comparing to \cite{b} we obtain here a non-trivial strictly singular operator only on a subspace of considered HI asymptotic $\ell_p$ spaces, nevertheless - thanks to the applied method - with much less restrictions on sets $N,L$ and parameters $(\theta_n)$, $(\rho_l)$ used in the construction of the spaces. 

Notice that the HI space in \cite{adkm} also admits a bounded strictly singular non-compact operators on a subspace by Th. 1.4 \cite{s2}, Prop. 3.3 \cite{adm} and the fact that its basis does not generate an $\ell_\omega$-spreading model. 
\er

\end{document}